\documentclass[12pt]{article}
\usepackage{epic,latexsym}
\usepackage{amssymb}
\usepackage{tikz}

\input epsf
\def \bbE {\mathbb{E}}

\newcommand{\diam}{{\rm diam}}

\newtheorem{thm}{Theorem}
\newtheorem{lem}[thm]{Lemma}
\newtheorem{cor}[thm]{Corollary}

\newtheorem{claim}{Claim}
\newtheorem{ob}{Observation}

\newtheorem{prob}{Problem}
\newtheorem{defn}{Definition}

\newcommand{\cF}{{\cal F}}

\newcommand{\cH}{{\cal H}}

\newcommand{\Vx}[1]{V(#1)}

\newcommand{\smallqed}{{\tiny ($\Box$)}}

\newcommand{\oo}{{\rm o}}

\newcommand{\gt}{\gamma_t}
\newcommand{\taut}{\tau_t}

\newcommand{\ONH}{{\rm ONH}}
\newcommand{\ec}{{\rm ec}}

\newcommand{\1}{\vspace{0.1cm}}
\newcommand{\2}{\vspace{0.2cm}}
\newcommand{\3}{\vspace{0.3cm}}

\newcommand{\proof}{\noindent\textbf{Proof.}\hspace{0.8em}}

\newcommand{\qed}{$\Box$}

\newenvironment{unnumbered}[1]{\trivlist \item [\hskip \labelsep {\bf
#1}]\ignorespaces\it}{\endtrivlist}

\def\vertex(#1){\put(#1){\circle*{2}}}
\def\vertexo(#1){\put(#1){\circle{2}}}
\def\vert(#1){\put(#1){\circle*{1.5}}}
\def\verto(#1){\put(#1){\circle{1.5}}}
\def\lab(#1)#2{\put(#1){\makebox(0,0)[c]{#2}}}
\def \nH {n_{_H}}
\def \mH {m_{_H}}

\def \nF {n_{_F}}

\newcommand{\TR}[1]{\mbox{$\tau(#1)$}}

\newcommand{\ST}[1]{\mbox{$\tau_s(#1)$}}

\def \nmr {\begin{enumerate}}
\def \enmr {\end{enumerate}}

\def \tmz {\begin{itemize}}
\def \etmz {\end{itemize}}

\begin{document}

\title{\vspace{-5ex} ~~ \\
 Total Transversals and Total Domination \\ in Uniform Hypergraphs}

\author{$\!\!\!\!\!\!\!\!\!\!^1$Csilla Bujt\'{a}s\thanks{Research supported in part by the Hungarian Scientific Research
  Fund, OTKA grant T-81493, and
 by the European Union and Hungary, co-financed
 by the European Social Fund through the project T\'AMOP-4.2.2.C-11/1/KONV-2012-0004 -- National Research Center
 for Development and Market Introduction of Advanced Information and Communication
 Technologies.}, \,
$^2$Michael A. Henning\thanks{Research supported in part by the
South African National Research Foundation and the University of
Johannesburg}, \, $^{1,3}$Zsolt Tuza$^*$ \, and $^{2,4}$Anders Yeo\\
\\
$^1$Department of Computer Science and Systems Technology \\
University of Pannonia \\
H-8200 Veszpr\'{e}m, Egyetem u.\ 10, Hungary\\
\small \tt Email: bujtas@dcs.vein.hu \\
\small \tt Email: tuza@dcs.vein.hu \\
\\
$^2$Department of Mathematics \\
University of Johannesburg \\
Auckland Park, 2006 South Africa\\
\small \tt Email: mahenning@uj.ac.za  \\
\\
$^3$~Alfr\'ed R\'enyi Institute of Mathematics \\
Hungarian Academy of Sciences \\
H-1053 Budapest, Re\'altanoda u.\ 13--15, Hungary \\
%\small \tt Email: \vspace{-3ex} tuza@sztaki.hu \\ \\
\\
$^4$Engineering Systems and Design \\
Singapore University of Technology and Design \\
20 Dover Drive Singapore, 138682, Singapore \\
\small \tt Email: \vspace{-3ex} andersyeo@gmail.com \\
~
}

\date{}
\maketitle

\newpage
\begin{abstract}
The first three authors [European J. Combin. 33 (2012), 62--71] established a relationship between the transversal number and the
domination number of uniform hypergraphs. In this paper, we establish a relationship between the total transversal number and the total domination number of uniform hypergraphs. We prove tight asymptotic upper bounds on the total transversal number in terms of the number of vertices, the number of edges, and the edge size.
\end{abstract}


\section{Introduction}

In this paper, we explore the study of total domination in
hypergraphs. We establish a relationship between the transversal
number and the total domination number of of uniform hypergraphs. We
introduce the concept of a total transversal in a hypergraph and
prove a general upper bound on the total domination number a uniform hypergraph in terms of its total transversal number.

Hypergraphs are systems of sets which are conceived as natural
extensions of graphs. A \emph{hypergraph} $H = (V(H),E(H))$ is a
finite set $V(H)$ of elements, called \emph{vertices}, together with
a finite multiset $E(H)$ of subsets of $V(H)$, called
\emph{hyperedges} or simply \emph{edges}. If the hypergraph $H$ is
clear from the context, we simply write $V = V(H)$ and $E = E(H)$. We
shall use the notation $\nH =|V|$ (or $n(H)$) and $\mH=|E|$ (or $m(H)$), and sometimes simply
$n$ and $m$ without subscript if the actual $H$ need not be
emphasized, to denote the \emph{order} and \emph{size} of $H$,
respectively. The edge set $E$ is often allowed to be a multiset in
the literature, but in the present context we exclude multiple edges.
Also, in the problems studied here, one may assume that $|\Vx{e}| \ge
2$ holds for all $e \in E$. An \emph{isolated edge} in $H$ is an edge
in $H$ that does not intersect any other edge in $H$. A \emph{linear
hypergraph} is a hypergraph in which every two edges intersect in at
most one vertex.

A $k$-\emph{edge} in $H$ is an edge of size~$k$. The hypergraph $H$
is said to be $k$-\emph{uniform} if every edge of $H$ is a $k$-edge.
The \emph{degree} of a vertex $v$ in $H$, denoted by $d_H(v)$ or
$d(v)$ if $H$ is clear from the context, is the number of edges of
$H$ which contain $v$. A vertex of degree $k$ is called a
\emph{degree-\,$k$ vertex}. The number of degree-$1$ vertices in $H$ is denoted by $n_1(H)$.
The minimum degree among the vertices of
$H$ is denoted by $\delta(H)$ and the maximum degree by $\Delta(H)$.

Two vertices $x$ and $y$ in $H$ are \emph{adjacent} if there is an edge $e$ of $H$ such that $\{x,y\} \subseteq \Vx{e}$. The \emph{neighborhood} of a vertex $v$ in $H$, denoted $N_H(v)$ or simply $N(v)$ if $H$ is clear from the context, is the set of all vertices different from $v$ that are adjacent to $v$. We call a vertex in $N(v)$ a \emph{neighbor} of $v$.
Two vertices $x$ and $y$ in $H$ are \emph{connected} if there is a
sequence $x=v_0,v_1,v_2\ldots,v_k=y$ of vertices of $H$ in which
$v_{i-1}$ is adjacent to $v_i$ for $i=1,2,\ldots,k$. A
\emph{connected hypergraph} is a hypergraph in which every pair of
vertices are connected. A maximal connected subhypergraph of $H$ is a
\emph{component} of $H$. Thus, no edge in $H$ contains vertices from
different components.

For a hypergraph $H$, the \emph{open neighborhood hypergraph} of $H$,
denoted by $\ONH(H)$, is the hypergraph with vertex set $V(H)$ and
edge set $\{N_H(v) \mid v \in V(H) \}$ consisting of the open
neighborhoods of vertices of $V(H)$ in $H$.

A subset $T$ of vertices in a hypergraph $H$ is a \emph{transversal} (also called \emph{vertex cover} or \emph{hitting set} in many papers) if $T$ has a nonempty intersection with every edge of $H$. The \emph{transversal number} $\TR{H}$ of $H$ is the minimum size of a transversal in $H$. A \emph{strong transversal}, often called a $2$-\emph{transversal}, in $H$ is a transversal that contains at least two vertices from every edge in $H$. The \emph{strong transversal number} $\ST{H}$ of $H$ is the minimum size of a strong transversal in $H$. Transversals in hypergraphs are well studied in the literature (see, for
example,~\cite{Al90,ChMc,DoHe13,HeLo12b,HeYe08,HeYe09,HeYe10,HeYe13,LaCh90,ThYe07,Tu90}).

We define a \emph{total transversal} in $H$ to be transversal $T$ in $H$ with the additional property that every vertex in $T$ has at least one neighbor in $T$, and we define the \emph{total transversal number} $\tau_t(H)$ of $H$ to be the minimum size of a total transversal in~$H$.

For a subset $X \subset V(H)$ of vertices in $H$, we define $H - X$ to be the hypergraph obtained from $H$ by deleting the vertices in $X$ and all edges incident with $X$, and deleting resulting isolated vertices, if any. We note that if $T'$ is a transversal in $H - X$, then $T' \cup X$ is a transversal in $H$. If $X = \{x\}$, then we write $H - X$ simply as $H - x$.

A \emph{dominating set} in a hypergraph $H = (V,E)$ is a subset of
vertices $D \subseteq V$ such that for every vertex $v \in V
\setminus D$ there exists an edge $e \in E$ for which $v \in {e}$ and ${e} \cap D \ne
\emptyset$. Equivalently, every vertex $v \in V \setminus D$ is
adjacent with a vertex in $D$. The \emph{domination number}
$\gamma(H)$ is the minimum cardinality of a dominating set in $H$. A
vertex $v$ in $H$ is said to be a \emph{dominating vertex} if it is
adjacent to every other vertex in $H$. A \emph{total dominating
set}, abbreviated TD-set, in a hypergraph $H = (V,E)$ is a subset of
vertices $D \subseteq V$ such that for every vertex $v \in V$ there
exists an edge $e \in E$ for which $v \in {e}$ and
 ${e} \cap (D \setminus \{v\}) \ne \emptyset$.
Equivalently, $D$ is a TD-set in $H$ if every vertex in $H$ is
adjacent with a vertex in $D$. The \emph{total domination number}
$\gamma_t(H)$ is the minimum cardinality of a TD-set in $H$. A
TD-set in $H$ of cardinality~$\gt(H)$ is called a $\gt(H)$-set.

While domination and total domination in graphs is very well studied in the literature (see, for example,~\cite{hhs1,hhs2,He09,HeYebook}), domination in hypergraphs was introduced relatively
recently by Acharya~\cite{Ac07} and studied further in
\cite{Ac08,BuHeTu12,HeLo12a,JoGeAb,JoTu09} and elsewhere.

A \emph{$2$-section graph}, $(H)_2$, of a hypergraph $H$ is defined as the graph with the same vertex set as $H$ and in which two edges are adjacent in $(H)_2$ if and only if they belong to a common edge in $H$.

Let $G$ be a graph. The \emph{degree} of a vertex $v$ in $G$ is denoted by $d_G(v)$ or
$d(v)$ if $G$ is clear from the context. The minimum degree among the vertices of
$G$ is denoted by $\delta(G)$ and the maximum degree by $\Delta(G)$.
An \emph{edge}-\emph{cover} in $G$ is a set of edges such
that every vertex in $G$ is incident with at least one edge in the edge-cover. We define a \emph{total edge}-\emph{cover} in $G$ to be an edge-cover that induces a subgraph with no isolated edge. We let $\ec_t(G)$ denote the minimum cardinality of a total edge-cover in
$G$. For two vertices $u$ and $v$ in a connected graph $G$, the
\emph{distance} $d(u,v)$ between $u$ and $v$ is the length of a
shortest $u$-$v$ path in $G$. The maximum distance among all pairs of vertices of $G$ is
the \emph{diameter} of $G$, which is denoted by $\diam(G)$.
A path and a cycle on $n$ vertices is denoted by $P_n$ and $C_n$, respectively.

The interplay between total domination in graphs and transversals in
hypergraphs has been studied in several papers (see, for example,
\cite{HeYe08,HeYe09,ThYe07}). The first three authors~\cite{BuHeTu12}
establish a relationship between the transversal number and the
domination number of uniform hypergraphs. In the present work, we establish a relationship between the total transversal number and the
total domination number of uniform hypergraphs. %We establish general upper bounds on the total domination number of a uniform hypergraph with minimum degree at least one in terms of its total transversal number, order and size.

\subsection{Key Definitions}
\label{S:defn}

We shall need the following definitions.

\begin{defn}
For an integer $k \ge 2$, let $\cH_k$ be the class of all $k$-uniform hypergraphs containing no isolated vertices or isolated edges
or multiple edges. Further, for $k \ge 3$ let $\cH_k^*$ consist of all hypergraphs in $\cH_k$ that have no two edges intersecting in $k-1$ vertices. We note that $\cH_k^*$ is a proper subclass of $\cH_k$.
 \label{defn1}
\end{defn}

\begin{defn}
For an integer $k \ge 2$, let
\[
b_k = \sup_{H \in \cH_k} \, \frac{ \taut(H) }{\nH + \mH}.
\]
 \label{defn2}
\end{defn}

\begin{defn}
For $k \ge 2$, let $H$ be obtained from a hypergraph $F \in \cH_k$ as follows. For each vertex $v$ in $F$, add $k$ new vertices $v_1,v_2,\ldots,v_k$ and two new $k$-edges $\{v,v_1,\ldots,v_{k-1}\}$ and $\{v_1,v_2,\ldots,v_k\}$.
%We note that these two added edges intersect in $k-1$ vertices.
Let $\cF_k$ denote the family of all such hypergraphs~$H$.
 \label{defn3}
\end{defn}

\begin{defn}
For $k \ge 3$, let $H$ be obtained from a hypergraph $F \in \cH_k^*$ as follows. For each vertex $v$ in $F$, add $k+1$ new vertices $v_1,v_2,\ldots,v_{k+1}$ and two new $k$-edges $\{v,v_1,v_2,\ldots,v_{k-1}\}$ and $\{v_2,v_3,\ldots,v_{k+1}\}$.
%We note that these two added edges intersect in $k-2$ vertices.
Let $\cF_k^*$ denote the family of all such hypergraphs $H$.
 \label{defn4}
\end{defn}

\section{Main Results}

We shall prove the following upper bounds on the total domination number of a uniform hypergraph in terms of its total transversal number, order and size. A proof of Theorem~\ref{t:main1A} is presented in Section~\ref{S:main1A}.

\begin{thm}
For $k \ge 3$, if $H \in \cH_k$, then $\displaystyle{\gt(H) \le \left( \max \left\{ \frac{2}{k+1}, b_{k-1}\right \} \right) \nH .}$
 \label{t:main1A}
\end{thm}

In view of Theorem~\ref{t:main1A}, it is of interest to determine the value of $b_k$ for $k \ge 2$. A proof of Theorem~\ref{t:main2a}
is presented in Section~\ref{S:main2a}.

\begin{thm}
$b_2 = \frac{2}{5}$, $b_3 = \frac{1}{3}$, and $b_4 \le \frac{1}{3}$. Further for $k \ge 5$, we have $b_k \le \frac{2}{7}$.
 \label{t:main2a}
\end{thm}

By Theorem~\ref{t:main2a}, we observe that
\[
b_{k-1} \le \frac{2}{k+1} \hspace*{1cm} \mbox{for $k \in \{3,4,5,6\}$.}
\]
Hence as a consequence of Theorem~\ref{t:main1A}  and Theorem~\ref{t:main2a}, and the well-known fact (see,~\cite{CoDaHe80}) that if $H \in \cH_2$, then $\gt(H) \le  2\nH/3$, we have the following result. The sharpness of the bound in Theorem~\ref{t:main2} is shown in Observation~\ref{ob3} in Section~\ref{S:known}.
%Recall that for $k \ge 2$, if $H \in \cF_k \subset \cH_k$, then $\gt(H) = 2\nH/(k+1)$, implying that the upper bound in Corollary~\ref{t:main2} is sharp.

\begin{thm}
For
%$2 \le k \le 6$,
$k \in \{2,3,4,5,6\}$,
if $H \in \cH_k$, then $\gt(H) \le 2\nH/(k+1)$, and this bound is sharp.
%$\displaystyle{ \gt(H) \le \left( \frac{2}{k+1} \right) \nH }$.
 \label{t:main2}
\end{thm}

The following result is a strengthening of the upper bound of Theorem~\ref{t:main1A} if we restrict the edges to intersect in at most~$k-2$ vertices. A proof of Theorem~\ref{t:main1B} is presented in Section~\ref{S:main1B}

\begin{thm}
For $k \ge 4$, if $H \in \cH_k^*$, then $\displaystyle{\gt(H) \le \left( \max \left\{ \frac{2}{k+2}, b_{k-1}\right \} \right) \nH .}$
 \label{t:main1B}
\end{thm}

\begin{cor}
For $k \ge 4$, if $H \in \cH_k^*$, then
$\gt(H) \le \nH/3$.
 \label{t:main3}
\end{cor}

The following result establishes a tight asymptotic  bound on $b_k$ for $k$ sufficiently large. A proof of Theorem~\ref{t:asym}
is presented in Section~\ref{S:asym}.

\begin{thm}
For $k$ sufficiently large, we have that
$\displaystyle{ b_k = (1 + \oo(1)) \frac{ \ln(k)}{k}.}$
 \label{t:asym}
\end{thm}

Theorem~\ref{t:asym} implies that the inequality $b_{k-1} \le 2/(k+1)$ is not valid when $k$ is large enough. This in turn, together with Theorem~\ref{t:main1A}, implies that Theorem~\ref{t:main2} is not true for large $k$.

\section{Known Results and Observations}
\label{S:known}

Cockayne et al.~\cite{CoDaHe80} established the following bound on
the total domination number of a connected graph in terms of its
order.

\begin{thm}{\rm (\cite{CoDaHe80})}
If $G$ is a connected graph of order $n \ge 3$, then $\gt(G) \le
2n/3$.
 \label{t:delta1}
\end{thm}

We shall need the following result due to Kelmans and  Mubayi~\cite{KeMu04}.

\begin{thm} {\rm (\cite{KeMu04})}
A cubic graph $G$ contains at least $\lceil |V(G)|/4 \rceil$ vertex disjoint $P_3$'s.
 \label{t:km}
\end{thm}

The following result shows that the total domination number of a
hypergraph $H$ is precisely the total domination of its $2$-section
graph and the transversal number of its open neighborhood hypergraph.

\begin{ob}
Let $H$ be a hypergraph with no isolated vertex. Then the following
holds. \\
\indent {\rm (a)} $\gt(H) = \gt((\cH)_2)$. \\
\indent {\rm (b)} $\gt(H) = \tau(\ONH(H))$. \label{ob1}
\end{ob}
\proof (a) Part~(a) follows readily from the fact that two vertices
in $H$ are adjacent in $H$ if and only if they are adjacent
in the $2$-section graph $(H)_2$ of $H$.

(b) On the one hand, every TD-set in $H$ contains a vertex from the
open neighborhood of each vertex in $H$, and is therefore a
transversal in $\ONH$, implying that $\tau(\ONH(H)) \le \gt(H)$. On
the other hand, every transversal in $\ONH$ contains a vertex from
the open neighborhood of each vertex of $H$, and is therefore a
TD-set in $G$, implying that $\gt(H) \le \tau(\ONH(H))$.
Consequently, $\gt(H) = \tau(\ONH(H))$.~\qed

\2
We shall need the following properties of hypergraphs in the family~$\cH_k$.

\begin{ob}
For $k \ge 2$, if $H \in \cH_k$, then the following hold. \\
\indent {\rm (a)} $\nH \ge k+1$, $\mH \ge 2$ and $\Delta(H) \ge 2$. \\
\indent {\rm (b)} $2\nH - n_1(H) \ge 2k$.
\label{ob2}
\end{ob}
\proof Part~(a) is immediate from the definition of the family~$\cH_k$. To prove Part~(b), let $n_{\ge 2}(H)$ denote the number of vertices in $H$ of degree at least~$2$. Let $e$ and $f$ be any two intersecting edges in $H$ and suppose they intersect in $\ell$ vertices, and so $|e \cup f| = 2k - \ell$. Then, $\nH \ge 2k - \ell \ge 2k - n_{\ge 2}(H)$, or, equivalently, $2\nH - n_1(H) = \nH + n_{\ge 2}(H) \ge 2k$.~\qed

\begin{ob}
The following holds. \\
\indent {\rm (a)} For $k \ge 2$, if $H \in \cF_k$, then $\gt(H) = 2\nH/(k+1)$. \\
\indent {\rm (b)} For $k \ge 3$, if $H \in \cF_k^*$, then $\gt(H) = 2\nH/(k+2)$.
\label{ob3}
\end{ob}
\proof For $k \ge 2$, let $H \in \cF_k$ be constructed as in Definition~\ref{defn3}. Then, $H \in \cH_k$ and $\nH = (k+1)\nF$. Every TD-set in $H$ contains at least two vertices in $\{v,v_1,v_2,\ldots,v_{k}\}$, implying that $\gt(H) \ge 2\nF$. However, the set $V(F) \cup T$, where $|T| = \nF$ and $T \subseteq V(H) \setminus V(F)$ consists of one added neighbor of each vertex in $V(F)$, is a TD-set in $H$, implying that $\gt(H) \le 2\nF$. Consequently, $\gt(H) = 2\nF = 2\nH/(k+1)$. For $k \ge 3$, let $H \in \cF_k^*$ be constructed as in Definition~\ref{defn4}. Then, $H \in \cH_k^*$ and $\gt(H) = 2\nF = 2\nH/(k+2)$.~\qed

%\newpage
\section{Preliminary Result}

We show first that total transversals of a $2$-regular hypergraph $H$ correspond to total
edge-covers in the dual multigraph, $G_H$, of $H$, where the
vertices of $G_H$ are the edges of $H$ and the edges of $G_H$
correspond to the vertices of $H$: if a vertex of $H$ is contained
in the edges $e$ and $f$ of $H$, then the corresponding edge of the
multigraph $G_H$ joins vertices $e$ and $f$ of $G_H$.

\begin{lem}
If $H$ is a linear $2$-regular hypergraph and $G_H$ is the dual of $H$, then $\taut(H) = \ec_t(G_H)$.
 \label{Dual}
\end{lem}
\proof By the linearity of $H$, the multigraph $G_H$ is in fact a
graph. Let $T$ be a total transversal in $H$ and let
$e$ be an arbitrary edge in $H$. Then there is a vertex $v \in T$
that covers~$e$. Since $H$ is $2$-regular, there is an edge $f$
different from $e$ that contains $v$. But then the edge in $G_H$
corresponding to the vertex $v$ in $H$ joins the two vertices $e$ and $f$ in $G_H$. Thus the edges of $G_H$ corresponding to vertices in $T$ form an
edge-cover in $G_H$.  Further suppose $u$ and $v$ are neighbors in
$H$ that belong to $T$ and let $g$ be the edge in $H$ containing $u$
and $v$. Let $e_u$ and $e_v$ be the edges, distinct from $g$, in $H$
containing $u$ and $v$. Then the edge in $G_H$ corresponding to the
vertex $u$ in $H$ joins the two vertices $e_u$ and $g$ in $G_H$,
while the edge in $G_H$ corresponding to the vertex $v$ in $H$ joins
the two vertices $e_v$ and $g$ in $G_H$, implying that the edges in
$G_H$ corresponding to the vertices $u$ and $v$ in $H$ have a vertex
in common. This implies that the edge-cover in $G_H$ corresponding to
the total transversal $T$ in $H$ is a total edge-cover in $G_H$.
Similarly, every total edge-cover in $G_H$ corresponds to a total
transversal in $H$. Therefore, $\taut(H) = \ec_t(G_H)$.~\qed

\newcommand{\DELETEtheBELOW}{
\medskip
Let $T_1$ denote the tree of order~$7$ obtained from a star $K_{1,3}$ by subdividing each edge exactly once, and let $T_2$ denote the tree of order~$7$ obtained from the disjoint union of two stars $K_{1,3}$ by identifying one leaf from each star. The trees $T_1$ and $T_2$ are shown in Figure~\ref{T1T2}(a) and~\ref{T1T2}(b), respectively.
%from three disjoint paths each on three vertices by identifying one end from each path. %Equivalently, $T_7^*$ is obtained from a star $K_{1,3}$ by subdividing each edge exactly %once.

\begin{figure}[htb]
\tikzstyle{every node}=[circle, draw, fill=black!0, inner sep=0pt,minimum width=.2cm]
\begin{center}
\begin{tikzpicture}[thick,scale=.7]
  \draw(0,0) { % <-- START CO-ORDINATES
    %%%%EDGES
    +(1.00,2.00) -- +(0.00,1.00)
    +(0.00,1.00) -- +(0.00,0.00)
    +(1.00,2.00) -- +(1.00,1.00)
    +(1.00,1.00) -- +(1.00,0.00)
    +(1.00,2.00) -- +(2.00,1.00)
    +(2.00,1.00) -- +(2.00,0.00)
    +(6.50,2.00) -- +(7.50,1.00)
    +(7.50,1.00) -- +(7.00,0.00)
    +(7.50,1.00) -- +(8.00,0.00)
    +(6.50,2.00) -- +(5.50,1.00)
    +(5.50,1.00) -- +(5.00,0.00)
    +(5.50,1.00) -- +(6.00,0.00)
    %%%%VERTICES
    +(1.00,2.00) node{}
    +(0.00,1.00) node{}
    +(1.00,1.00) node{}
    +(2.00,1.00) node{}
    +(2.00,0.00) node{}
    +(1.00,0.00) node{}
    +(0.00,0.00) node{}
    +(5.50,1.00) node{}
    +(6.00,0.00) node{}
    +(5.00,0.00) node{}
    +(8.00,0.00) node{}
    +(7.00,0.00) node{}
    +(7.50,1.00) node{}
    +(6.50,2.00) node{}
    %%%%LABELS
    +(1.20,-0.75) node[rectangle, draw=white!0, fill=white!100]{(a) $T_1$}
    +(6.50,-0.75) node[rectangle, draw=white!0, fill=white!100]{(b) $T_2$}
     %%% <-- BELOW
    %+(-0.75,1.00) node[rectangle, draw=white!0, fill=white!100]{$G$}   %%% <-- LEFT
    %+(4.00,2.75) node[rectangle, draw=white!0, fill=white!100]{$G$}   %%% <-- ABOVE
    %+(8.75,1.00) node[rectangle, draw=white!0, fill=white!100]{$G$}   %%% <-- RIGHT
  };
\end{tikzpicture}
\end{center}
\vskip -0.6 cm \caption{The trees $T_1$ and $T_2$.} \label{T1T2}
\end{figure}

\begin{lem}
Let $T$ be a tree of order at least~$7$ and with $\Delta(T) \le 3$.
If $T \ne T_1$, then there exists an edge $e$ of $T$ such that both
components of $T - e$ have order at least~$3$. \label{l:tree}
\end{lem}
\proof  Let $T \ne T_1$ be a tree of order~$n \ge 7$ and with $\Delta(T) \le 3$. If $\diam(T) \ge 5$, then let $e$ be a middle edge on a longest path in $T$. Then both components of $T - e$ have order at least~$3$. Hence we may assume that $\diam(T) \le 4$, for otherwise the desired result follows. In particular, we note that $T \ne P_n$, implying that $\Delta(T) = 3$. Let $v$ be a vertex of degree~$3$ in $T$. If $v$ has a neighbor $u$ of degree~$3$, then taking $e = uv$ the desired result follows. Hence we may assume that every neighbor of $v$ has degree at most~$2$ in $T$. If every neighbor of $v$ has degree~$2$, then the constraint that the diameter is at most~$4$ implies that $T = T_1$, a contradiction. If $v$ has only one neighbor that is a leaf, then the constraint that the order is at least~$7$ implies that there is a path in $T$ of length at least~$5$, a contradiction. Hence exactly two neighbors of $v$ are leaves. The constraints that the diameter is at most~$4$ and the order at least~$7$ imply that $T = T_2$. But then if $e$ denotes one of the two edges incident with the central vertex of $T$, then both components of $T - e$ have order at least~$3$.~\qed
}

\section{Proof of Main Results}

\subsection{Proof of Theorem~\ref{t:main1A}}
\label{S:main1A}

In this section, we present a proof of Theorem~\ref{t:main1A}. Recall
its statement.

\noindent \textbf{Theorem~\ref{t:main1A}}. \emph{For $k \ge 3$, if $H \in \cH_k$, then $\displaystyle{\gt(H) \le \left( \max \left\{ \frac{2}{k+1}, b_{k-1}\right \} \right) \nH .}$
}

\noindent \textbf{Proof of Theorem~\ref{t:main1A}.} Suppose to the contrary that the theorem is not true. Let $H \in \cH_k$ be a counterexample with $\nH + \mH$ a minimum. In what follows we present a series of claims describing some structural properties of $H$ which culminate in the implication of its non-existence.

\begin{claim}\label{claim1}
The following properties hold in the hypergraph $H$. \\
\indent {\rm (a)} $H$ is connected. \\
\indent {\rm (b)} The deletion of any edge in $H$ creates an isolated vertex or an isolated edge. \\
\indent {\rm (c)} There is no dominating vertex in $H$.
\end{claim}
\textbf{Proof of Claim~\ref{claim1}.} Part~(a) is immediate from the
minimality of $H$. Part~(b) is also immediate since the deletion of
an edge cannot decrease the total domination number. To prove
Part~(c), suppose that $H$ contains a dominating vertex $v$. The
vertex $v$ and any one of its neighbors forms a TD-set in $H$,
implying that $\gt(H) = 2$. As $H \in \cH_k$, there is no
isolated vertex or isolated edge in $H$, implying that $\nH \ge
k+1$. Hence, $\gt(H) \le 2\nH/(k+1)$, contradicting the fact that
$H$ is a counterexample to the theorem. This proves
Part~(c).~\smallqed

\begin{claim}\label{claim2}
Every edge in $H$ contains at least one degree-$1$ vertex.
\end{claim}
\textbf{Proof of Claim~\ref{claim2}.} Suppose to the contrary that
there is an edge $e$ that does not contain any degree-$1$ vertices.
Thus every vertex contained in $e$ has degree at least~$2$ in $H$.
By Claim~\ref{claim1}(b), there is therefore an edge, $e_1$, which
would become isolated after the deletion of the edge $e$ from
$H$. Thus, every vertex in $e \cap e_1$ has
degree~$2$ in $H$, while every vertex in $e_1 \setminus e$ has
degree~$1$ in $H$. Let $v \in e \cap e_1$. Then, $d_H(v) = 2$. By
Claim~\ref{claim1}(a), $H$ is connected and by
Claim~\ref{claim1}(c), the vertex $v$ is not a dominating vertex of
$H$, implying that there exists an edge, $e_2$, such that $v \notin
e_2$ but $e_2$ intersects $e$. Since $e \ne e_2$ and $v \notin e_2$,
we note that $e_1 \cap e_2 = \emptyset$. Let $u \in e \cap e_2$ and
note that $u \notin e_1$.

Initially we set $T = \emptyset$ and we now construct a hypergraph
$H'$ from $H$ as follows. We delete all edges incident with $u$ or
$v$ or with both $u$ and $v$ and we delete any resulting isolated
vertices. Further we add both vertices $u$ and $v$ to the set $T$.
We note that the edges $e$ and $e_1$ are both deleted, implying that
every vertex in $e_1$ becomes an isolated vertex. Further since we
remove all edges incident with $u$, the vertex $u$ becomes an
isolated vertex. We therefore delete at least $k+1$ vertices and we
add two vertices to $T$. If this process creates an isolated edge,
then such an isolated edge necessarily contains a vertex that is
adjacent to at least one of $u$ and $v$ (for otherwise it would be
an isolated edge in $H$, a contradiction).
From each such isolated edge $f$, if any, we choose one vertex that
is a neighbor of $u$ or $v$ and add it into $T$, and delete the $k$
vertices in $f$. Hence, $|T| = 2 + \ell$, where $\ell \ge 0$ denotes
the number of isolated edges created when removing $u$ and $v$.

Let $n'$ denote the number of vertices in $H$ that are not deleted in the process (possibly, $n' = 0$). At least $k + 1 + k \ell$ vertices were deleted from $H$. Thus, $n' \le \nH - k - 1 - k \ell$, implying that

\[
\begin{array}{lcl} \3
\displaystyle{ \left( \frac{2}{k+1} \right) (\nH - n') } & \ge &
 \displaystyle{ \left( \frac{2}{k+1} \right) \left( k+1 + k \ell  \right) } \\ \3
& = & \displaystyle{  2 +  \left( \frac{2k}{k+1} \right) \ell }  \\ \3
& \ge & \displaystyle{  2 + \ell } \\ \3
& = & |T|.
\end{array}
\]

If $n' = 0$, then the set $T$ is a TD-set in $H$, implying that $\gt(H) \le |T| \le 2\nH/(k+1)$, a contradiction. Hence, $n' > 0$. Let $H'$ denote the resulting hypergraph on these $n'$ vertices. Let $H'$ have size~$m'$. By construction, the hypergraph $H'$ is in the family~$\cH_k$. In particular, we note that $n' \ge k+1$. By the minimality of $H$, we have that
\[
\gt(H') \le \left( \max \left\{ \frac{2}{k+1}, b_{k-1}\right \} \right) n'.
\]

Let $T'$ be a $\gt(H')$-set and note that the set $T \cup T'$ is a TD-set of $H$. Suppose that $2/(k+1) \ge b_{k-1}$. Then, $|T'| \le 2n'/(k+1)$, and so
\[
\gt(H) \le |T \cup T'| \le \left( \frac{2}{k+1} \right) (\nH - n') + \left( \frac{2}{k+1} \right) n' = \left( \frac{2}{k+1} \right) \nH,
\]

%\noindent
a contradiction. Hence, $2/(k+1) < b_{k-1}$. Thus, $|T'| \le b_{k-1}n'$, and so
\[
\gt(H) \le |T \cup T'| \le \left( \frac{2}{k+1} \right) (\nH - n') + b_{k-1}n' < b_{k-1} (\nH - n') + b_{k-1}n'  = b_{k-1} \nH,
\]

a contradiction. This completes the proof of Claim~\ref{claim2}.~\smallqed

\2
We now return to the proof of Theorem~\ref{t:main1A}. By Claim~\ref{claim2}, every edge in $H$ contains at least one degree-$1$ vertex. If there are two edges, $f_1$ and $f_2$, in $H$ that intersect in $k-1$ vertices, then for $j \in \{1,2\}$, the edge $f_{j}$ contains exactly one vertex, $v_j$ say, not in $f_{3-j}$ and this vertex has degree~$1$ in $H$. Thus if we delete the vertices $v_1$ and $v_2$ from $H$, then we would create a multiple edge, namely $f_1' = f_1 \setminus \{v_1\}$ and $f_2' = f_2 \setminus \{v_2\}$.
Let $H'$ be the hypergraph obtained from $H$ by deleting exactly one degree-1 vertex from each edge and by replacing resulting multiple edges, if any, by single edges. Let $H'$ have order $n'$ and size $m'$. Then, $n' = \nH - \mH$ and $m' \le  \mH$. Thus, $n' + m' \le \nH$.

\begin{claim}\label{claim3}
$H' \in \cH_{k-1}$ and $\taut(H') \le b_{k-1} \nH$.
\end{claim}
\textbf{Proof of Claim~\ref{claim3}.} If $H'$ contains an isolated edge, then every vertex in such an isolated edge would be
a dominating vertex in $H$, contradicting Claim~\ref{claim1}(c). Hence, $H'$ contains no isolated edge. By construction, $H'$ has no multiple edges and no isolated vertices. Therefore, $H' \in \cH_{k-1}$. We note that $k - 1 \ge 2$. By Definition~\ref{defn2} we have that
$\taut(H') \le (n' + m') b_{k-1} \le b_{k-1} \nH$.~\smallqed

\begin{claim}\label{claim4}
$\taut(H')=\gt(H)$.
\end{claim}
\textbf{Proof of Claim~\ref{claim4}.}  Among all $\gt(H)$-sets, let $S$ be chosen to contain as few vertices of degree~$1$ in $H$ as possible. Suppose that $S$ contains a degree-1 vertex, $x$, in $H$. Let $e_x$ be the edge containing~$x$. By the minimality of the set $S$, the set $S_x = S \setminus \{x\}$ is not a TD-set in $H$. Let $y$ be a vertex in $S$ that is adjacent to $x$ in $H$. Then, $y \in e_x$. If $y$ is adjacent to a vertex of $S_x$, then the set $S_x$ would be a TD-set in $H$, a contradiction. Hence, $y$ is adjacent to no vertex of $S$ except for the vertex $x$. Since $H$ contains no dominating vertex and since $H$ has no isolated edge, there exists a neighbor, $w$ say, of $y$ that has degree at least~$2$ in $H$. But then $S_x \cup \{w\}$ is a TD-set of $H$ of cardinality~$|S| = \gt(H)$ that contains fewer degree-$1$ vertices than does $S$, contradicting our choice of the set $S$. Therefore, $S$ contains no vertices of degree~$1$, implying that $S \subseteq V(H')$. Further if $S$ is not a transversal in $H$, then let $e'$ be an edge in $H$ not intersected by $S$. But since $e'$ contains a degree-1 vertex, such a vertex would not be (totally) dominated by $S$ in $H$, a contradiction. Hence, $S$ is a transversal in $H$. Further since every vertex in the TD-set $S$ has a neighbor in $H$ that belongs to $S$, the set $S$ is in fact a total transversal of $H$. Since $S \subseteq V(H')$, the set $S$ is therefore also a total transversal of $H'$, implying that $\taut(H') \le \gt(H)$. Conversely, every total transversal in $H'$ is a TD-set in $H'$ and therefore also in $H$, implying that $\gt(H) \le \taut(H')$. Consequently, $\taut(H')=\gt(H)$.~\smallqed

\2
By Claim~\ref{claim3} and Claim~\ref{claim4}, we have that $\gt(H) \le b_{k-1} \nH$, a contradiction. This  completes the proof of Theorem~\ref{t:main1A}.~\qed

\subsection{Proof of Theorem~\ref{t:main2a}}
\label{S:main2a}

In this section, we present a proof of Theorem~\ref{t:main2a}. We first consider the family~$\cH_2$.

\begin{thm}
%$b_2 = 2/5$.
If $H \in \cH_2$, then $\taut(H) \le 2(\nH+ \mH)/5$.
 \label{t:b2}
\end{thm}
\textbf{Proof of Theorem~\ref{t:b2}.} Suppose to the contrary that the theorem is not true.
Let $H \in \cH_2$ be a counterexample with $\nH + \mH$ a minimum. Clearly, $H$ is connected.
By Observation~\ref{ob2}, we have that $\nH \ge 3$, $\mH \ge 2$ and $\Delta(H) \ge 2$. If $\taut(H) = 2$, then the result is immediate. Hence we may assume that $\taut(H) \ge 3$. Let $x$ be a vertex of maximum degree in $H$. Since $\taut(H) \ge 3$, there is a neighbor $y$ of $x$ that is not isolated in $H - x$. We delete the vertices $x$ and $y$ and all edges incident with $x$ or $y$, together with any resulting isolated vertices, if any, and let $T = \{x,y\}$. Further if this process creates an isolated edge, $e$, then such an isolated edge necessarily contains a vertex that is adjacent to $x$ {or} $y$, for otherwise the edge $e$ would be an isolated edge in $H$, a contradiction. From each such isolated edge $e$, if any, we choose one vertex that is a neighbor of $x$ or $y$ and add it to the set $T$, and delete the two vertices in $e$. Suppose that $\ell \ge 0$ isolated edges were created when $x$ and $y$ are deleted. Then, $|T| = 2 + \ell$ and at least $2 + 2\ell$ vertices and at least $3 + \ell$ edges were deleted. Let $H'$ denote the resulting graph. Thus, if $H'$ has $n'$ vertices and $m'$ edges, then $n' + m' \le \nH + \mH - (5 + 3\ell)$. Since $H$ is a minimum counterexample, we have that $\taut(H') \le 2(n' + m')/5$, implying that

\[
\begin{array}{lcl} \1
\taut(H)  & \le &\taut(H') + |T|  \\ \1
& \le & \frac{2}{5}(\nH + \mH - 5 - 3\ell) + 2 + \ell \\ \1
& \le & \frac{2}{5}(\nH + \mH) - \frac{\ell}{5} \\ \1
& \le & \frac{2}{5}(\nH + \mH),
\end{array}
\]
contradicting the fact that $H$ is a counterexample.~\qed

\2
As an immediate consequence of Theorem~\ref{t:b2}, we have that $b_2 \le 2/5$.
Taking $H$ to be a path $P_3$ on three vertices, we note that $H \in \cH_2$ and $\taut(H) = 2 = 2(\nH+ \mH)/5$, implying that $b_2 \ge 2/5$. Consequently, $b_2 = 2/5$. This can also be seen by considering the cycle of order five, $C_5$, instead of $P_3$, as $\taut(C_5)=4$. We state this formally as follows.

\begin{cor}
$b_2 = 2/5$.
 \label{c:b2}
\end{cor}

We next consider the family~$\cH_k$, where $k \ge 3$.

\begin{thm}
For $k \ge 3$, if $H \in \cH_k$, then $\taut(H) \le (\nH+ \mH)/3$.
 \label{t:kge3}
\end{thm}
\textbf{Proof of Theorem~\ref{t:kge3}.} Suppose to the contrary that the theorem is not true. Let $H \in \cH_k$ be a counterexample with $\nH + \mH$ a minimum. Clearly, $H$ is connected since otherwise the theorem holds for each component of $H$ and therefore also for $H$, a contradiction. By Observation~\ref{ob2}, we have that $\nH \ge k+1$, $\mH \ge 2$ and $\Delta(H) \ge 2$. In what follows we present a series of claims describing some structural properties of $H$ which culminate in the implication of its non-existence.

\begin{unnumbered}{Claim~A.}
$\taut(H) \ge 3$ and no vertex is incident with every edge in $H$.
\end{unnumbered}
\textbf{Proof of Claim~A.} Suppose to the contrary that $\taut(H) < 3$.
Then, $\taut(H) = 2$. Since $\nH + \mH \ge k+3 \ge 6$, we therefore have that $\taut(H) = 2 \le (\nH+\mH)/3$, contradicting the fact that $H$ is a counterexample. Hence, $\taut(H) \ge 3$.

If there is a vertex $v$ incident with every edge in $H$, then the vertex $v$ and one of its neighbors form a total transversal in $H$, implying that $\taut(H) = 2$, a contradiction.
Hence, no vertex is incident with every edge in $H$.~\smallqed

\begin{unnumbered}{Claim~B.}
There is no set $X \subset V(H)$, such that (a) and (b) below  hold. \\
      \indent {\rm (a)}  Every vertex in $X$ has a neighbor in $H$ in the set $X$.  \\
      \indent {\rm (b)} $n(H-X) + m(H-X) \le \nH + \mH - 3|X|$.
\end{unnumbered}
\textbf{Proof of Claim~B.} Suppose to the contrary that a subset $X \subset V(H)$ satisfying the two conditions in the statement of the claim exists. Let $H' = H - X$. By supposition, $n(H') + m(H') \le \nH + \mH - 3|X|$.

Let $e_1,\ldots,e_\ell$, where $\ell \ge 0$, be the isolated edges in $H'$. Since $H$ contains no isolated edge, each isolated edge in $H'$ contains a vertex of degree at least~$2$ in $H$. For each $i = 1,\ldots,\ell$, let $z_i \in e_i$ be chosen so that $d_H(z_i) \ge 2$, and let $X^* = X \cup \{z_1,\ldots,z_\ell\}$. We note that every vertex in $X^*$ is adjacent to a vertex in $X \subseteq X^*$.

Let $H^* = H - X^*$. By construction, $H^* \in \cH_k$. Moreover, $n(H^*) = n(H') - k \ell$ and $m(H^*) = m(H') - \ell$. By the minimality of $H$, we have that $\taut(H^*) \le (n(H^*)+m(H^*))/3$. Since every $\taut(H^*)$-set can be extended to a total transversal of $H$ by adding to it the set $X^*$, and since $k \ge 3$, we have that

\[
\begin{array}{lcl} \2
\taut(H) & \le & \taut(H^*) + |X^*| \\ \2
& \le & \displaystyle{ \frac{1}{3}  (n(H^*) + m(H^*)) + |X| + \ell } \\ \2
& {=} & \displaystyle{ \frac{1}{3}  (n(H') - k \ell + m(H') - \ell ) + |X| + \ell } \\ \2
& \le & \displaystyle{ \frac{1}{3}  (\nH+\mH -3|X|  - k \ell - \ell ) + |X| + \ell } \\ \2
& \le & \displaystyle{ \frac{1}{3}(\nH + \mH) }, \\
\end{array}
\]

\noindent
{contradicting the fact that $H$ is a counterexample}.~\smallqed

\begin{unnumbered}{Claim~C.}
$\Delta(H) = 2$.
\end{unnumbered}
\textbf{Proof of Claim~C.} Suppose to the contrary that $\Delta(H) \ge 3$. Let $x$ be a vertex of maximum degree in $H$. By Claim~A, the vertex $x$ is not incident with every edge in $H$. Hence since $H$ is connected, there exists an edge, $e$, that contains a neighbor, $y$, of $x$ but does not contain~$x$. Let $X=\{x,y\}$ and note that $n(H-X) \le \nH-2$ and $m(H-X) \le \mH - 4$. As $x$ and $y$ are adjacent in $H$, we obtain a contradiction to Claim~B.~\smallqed

\begin{unnumbered}{Claim~D.}
%$H$ is a $2$-regular $k$-uniform hypergraph
$H$ is $2$-regular.
\end{unnumbered}
\textbf{Proof of Claim~D.} Suppose that there exists a vertex $v_1$
of degree~$1$ in $H$. Let $e_1$ be the edge incident with $v_1$.
Since $H$ has no isolated edge, let $e_2$ be an edge intersecting
$e_1$, and let $v_2 \in e_1 \cap e_2$. By Claim~A, the vertex $v_2$ is not incident with every edge in $H$. Hence there exists an
edge, $e_3$, not containing $v_2$ that intersects $e_1$ or $e_2$ in a vertex $v_3$. Let $X = \{v_2,v_3\}$ and note that the vertices $v_1,v_2,v_3$ and the edges $e_1,e_2,e_3$ are removed from $H$ in order to create $H-X$. Therefore, $n(H-X) \le \nH-3$ and $m(H-X) \le \mH - 3$, which as $v_2$ and $v_3$ are adjacent in $H$, contradicts Claim~B.~\smallqed

\begin{unnumbered}{Claim~E.}
$H$ is a linear hypergraph.
\end{unnumbered}
\textbf{Proof of Claim~E.} By Claim~D, $H$ is a $2$-regular $k$-uniform hypergraph. Suppose that there are two edges $e$ and $f$ having two or more vertices in common. Let $v$ be a vertex in $e$ that does not belong to $e \cap f$. Since $H$ is $2$-regular, there is an edge $g$ which contains $v$ but is different from $e$ or $f$. Let $u$ be a vertex in $e \cap f$. Since $u$ and $v$ belong to the common edge $e$, they are neighbors in $H$. Let $X = \{u,v\}$ and note that the vertices in  $\{v\} \cup (e \cap f)$ and the edges $e,f,g$ are removed from $H$ in order to create $H-X$. Therefore, $n(H-X) \le \nH-3$ and $m(H-X) \le \mH - 3$, which contradicts Claim~B.~\smallqed

\2
By Claim~D and Claim~E, $H$ is a $2$-regular $k$-uniform linear connected hypergraph.

\begin{unnumbered}{Claim~F.}
$k=3$
\end{unnumbered}
\textbf{Proof of Claim~F.} Suppose to the contrary that $k \ge 4$.
Then, $\nH = k\mH/2 \ge 2\mH$. We now consider the dual, $G_H$, of the hypergraph $H$. By the 2-regularity and the linearity of $H$, the dual $G_H$ is a graph. Since $H$ is $k$-uniform, the graph $G_H$ is $k$-regular. Further since $H$ is
connected, so too is $G_H$. By construction, $G_H$ has order~$n(G_H) = \mH$ and size~$m(G_H) = \nH$. Let $T$ be a spanning tree in $G_H$.
Since the set $E(T)$ of edges of $T$ form a total edge-cover in $G_H$ and since $\nH \ge 2\mH$, we have by Lemma~\ref{Dual} that
$\taut(H) = \ec_t(G_H) \le |E(T)| = n(G_H) - 1 = \mH - 1 < \frac{1}{3} (\nH + \mH)$, a contradiction.~\smallqed

\2
By Claim~D,~E and~F, we have that $H$ is a $2$-regular $3$-uniform linear connected hypergraph. We now consider the dual, $G_H$, of the hypergraph $H$. We note that the dual, $G_H$, is a connected, cubic graph. Applying Theorem~\ref{t:km} to the cubic graph $G_H$, there exist at least  $\lceil n(G_H)/4 \rceil$ vertex disjoint $P_3$'s in $G_H$. Let $G_1,G_2,\ldots,G_{\ell}$ denote vertex disjoint subgraphs in $G_H$ each of which are isomorphic to $P_3$, such that $\ell \ge \lceil n(G_H)/4 \rceil \ge m_H/4$. If some vertex does not belong to one of these subgraphs $G_1,G_2,\ldots,G_{\ell}$, then the connectivity of $G_H$ implies that there is an edge, $e$, joining a vertex in $V(G_i)$ for some $i$, $1 \le i \le \ell$, and a vertex, $x$, not belonging to any subgraph $G_1,G_2,\ldots,G_{\ell}$. We now add the vertex $x$ and edge $e$ to the subgraph $G_i$. We continue this process until all vertices in $G_H$ belong to exactly one of the resulting subgraphs $G_1,G_2,\ldots,G_{\ell}$. The subgraph of $G_H$ induced by the edges in these $\ell$ subgraphs is a spanning forest, $F$, of $G_H$, that contains~$\ell \ge \mH/4$ components each of which contain at least three vertices.

Since every component of $F$ has order at least~$3$, the set $E(F)$ of edges of $F$ forms a total edge-cover in $G_H$. Since $n(G_H) = \mH$ and $\ell \ge m_H/4$, we have that $|E(F)| = n(G_H) - \ell \le 3\mH/4$. Therefore, recalling that $\nH = 3\mH/2$, we have by Lemma~\ref{Dual} that

\[
\taut(H) = \ec_t(G_H) \le |E(F)| \le \frac{3}{4} \, \mH \le \frac{1}{3} (\nH + \mH),
\]

\noindent {a contradiction}. This completes the proof of Theorem~\ref{t:kge3}.~\qed

\2 As an immediate consequence of Theorem~\ref{t:kge3}, we have that
$b_k \le 1/3$ for all $k \ge 3$. Taking $H$ to be the hypergraph of
order~$\nH = 4$ and size~$\mH = 2$ where the two edges of $H$
intersect in two vertices, we note that $H \in \cH_3$ and
 $\taut(H) = 2 = (\nH+
\mH)/6$, implying that $b_3 \ge 1/3$. Consequently, $b_3 = 1/3$. As
observed earlier, $b_4 \le 1/3$. We state this formally as follows.

\begin{cor}
$b_3 = \frac{1}{3}$ and $b_4 \le \frac{1}{3}$.
 \label{c:b3b4}
\end{cor}

We remark that the result of Theorem~\ref{t:kge3} can be strengthened slightly when $k \ge 4$, as the following result shows. We omit the proof (which is similar, but simpler, to the proof of Theorem~\ref{t:kge5} presented below).

\begin{thm}
For $k \ge 4$, if $H \in \cH_k$, then $6\taut(H) \le 2\nH + 2\mH - n_1(H)$.
 \label{t:kge4}
\end{thm}

We next consider the family~$\cH_k$, where $k \ge 5$.

\begin{thm}
For $k \ge 5$, if $H \in \cH_k$, then $7\taut(H) \le 2\nH + 2\mH - n_1(H)$.
 \label{t:kge5}
\end{thm}
\textbf{Proof of Theorem~\ref{t:kge5}.} For $k \ge 5$ and all hypergraphs $H \in \cH_k$, let
\[
\Theta(H) = 2\nH + 2\mH - n_1(H).
\]

We wish to show that $7\taut(H) \le \Theta(H)$. Suppose to the contrary that the theorem is not true. Let $H \in \cH_k$ be a counterexample with minimum $\Theta(H)$. Clearly, $H$ is connected since otherwise the theorem holds for each component of $H$ and therefore also for $H$, a contradiction. By Observation~\ref{ob2}(a), we have that $\nH \ge k+1$, $\mH \ge 2$ and $\Delta(H) \ge 2$. By Observation~\ref{ob2}(b), we have that $2\nH - n_1(H) \ge 2k$. In what follows we present a series of claims describing some structural properties of $H$ which culminate in the implication of its non-existence.

\begin{unnumbered}{Claim~I.}
$\taut(H) \ge 3$.
\end{unnumbered}
\textbf{Proof of Claim~I.} Suppose that $\taut(H) < 3$. Then, $\taut(H) = 2$. Since $2\nH - n_1(H) \ge 2k$ and $\mH \ge 2$, we therefore have that $7\taut(H) = 14 \le 2k + 4 \le \Theta(H)$, {contradicting the fact that $H$ is a counterexample}.~\smallqed

\begin{unnumbered}{Claim~II.}
If $X$ is a set of vertices in $H$, such that every vertex in $X$ is adjacent to some other vertex of $X$, then $\Theta(H - X) > \Theta(H) - 7|X|$.
\end{unnumbered}
\textbf{Proof of Claim~II.} Suppose to the contrary that exists a subset $X \subset V(H)$ such that every vertex in $X$ is adjacent to some other vertex of $X$ but $\Theta(H - X) \le \Theta(H) - 7|X|$.
Let $H' = H - X$. Let $e_1,\ldots,e_\ell$, where $\ell \ge 0$, be the isolated edges in $H'$. Since $H$ contains no isolated edge, every isolated edge in $H'$ contains a vertex of degree at least~$2$ in $H$ that is adjacent to a vertex of $X$ in $H$. For each $i = 1,\ldots,\ell$, let $z_i \in e_i$ be chosen so that $d_H(z_i) \ge 2$, and let $X^* = \{z_1,\ldots,z_\ell\}$. We note that every vertex in $X \cup X^*$ is adjacent to some other vertex of $X$. We now consider the hypergraph $H^* = H' - X^*$.

We note that $H^* \in \cH_k$. When constructing $H^*$ from $H'$ we deleted all $k \ell$ vertices from the $\ell$ isolated edges in $H'$ and we deleted all $\ell$ isolated edges. Since each such deleted vertex has degree~$1$ in $H'$, the contribution of the $k \ell$ deleted vertices from $H'$ to the sum $2n(H') - n_1(H')$ is~$k \ell$. The contribution of the $\ell$ deleted edges to the sum $2m(H')$ is~$2\ell$. By supposition, $\Theta(H') \le \Theta(H) - 7|X|$. Since $k \ge 5$, we therefore have that
\[
\begin{array}{lcl} \1
\Theta(H^*) & = & \Theta(H') - \ell(k+2) \\ \1
& \le & \Theta(H') - 7\ell \\ \1
& \le & (\Theta(H) - 7|X|) - 7\ell \\
& = & \Theta(H) - 7|X| - 7|X^*|. \\
\end{array}
\]

By the minimality of $\Theta(H)$, we have that $7\taut(H^*) \le \Theta(H^*)$. Every (minimum) total transversal in $H^*$ can be extended to a total transversal in $H$ by adding to the set $X \cup X^*$, implying that $\taut(H) \le \taut(H^*) + |X| + |X^*|$. Hence,
\[
\begin{array}{lcl} \1
7\taut(H) & \le & 7\taut(H^*) + 7|X| + 7|X^*| \\ \1
& \le & \Theta(H^*) + 7|X| + 7|X^*| \\ \1
& \le & \Theta(H),
\end{array}
\]

\noindent
{a contradiction}.~\smallqed

\begin{unnumbered}{Claim~III.}
$\Delta(H) \le 3$.
\end{unnumbered}
\textbf{Proof of Claim~III.} Suppose to the contrary that $\Delta(H) \ge 4$. Let $x$ be a vertex of maximum degree in $H$. Since $\taut(H) \ge 3$ by Claim~I, and since $H$ is connected, there exists an edge, $e$, that contains a neighbor, $y$, of $x$ but does not contain~$x$. Let $X = \{x,y\}$ and consider the hypergraph $H - X$. Since $d_H(x) \ge 4$ and $d_H(y) \ge 2$, the vertices $x$ and $y$ both contribute~$2$ to the sum $2n(H) - n_1(H)$. Further since at least five distinct edges are deleted from $H$ when constructing $H - X$, the contribution of the deleted edges to the sum $2m(H)$ is at least~$10$. Hence, $\Theta(H - X) \le \Theta(H)  - 14 = \Theta(H)  - 7|X|$, contradicting Claim~II.~\smallqed

\begin{unnumbered}{Claim~IV.}
$\Delta(H) = 2$.
\end{unnumbered}
\textbf{Proof of Claim~IV.}
As observed earlier, $\Delta(H) \ge 2$. By Claim~III, $\Delta(H) \le 3$. Suppose to the contrary that $\Delta(H) = 3$. Let $x$ be a vertex with $d_H(x) = 3$ and consider the hypergraph $H' = H-x$.
Suppose that $d_{H'}(y) \ge 2$ for some $y \in N_H(x)$.
Let $X = \{x,y\}$ and consider the hypergraph $H - X$.
Since $d_H(x) = 3$ and $d_H(y) = 3$, the vertices $x$ and $y$ both contribute~$2$ to the sum $2n(H) - n_1(H)$. Further since five distinct edges are deleted from $H$ when constructing $H - X$, the contribution of the deleted edges to the sum $2m(H)$ is~$10$. Hence, $\Theta(H - X) \le \Theta(H)  - 14 = \Theta(H) - 7|X|$, contradicting Claim~II. Therefore, $d_{H'}(y) \le 1$ for every vertex $y \in N_H(x)$.

Since $\taut(H) \ge 3$ by Claim~I, and since $H$ is connected, there exists a neighbor, $y^*$, of $x$ that has degree at least~$1$ in $H'$. Let $X^* = \{x,y^*\}$ and consider the hypergraph $H^* = H - X^*$. Since $d_H(x) = 3$ and $d_H(y^*) \ge 2$, the vertices $x$ and $y^*$ both contribute~$2$ to the sum $2n(H) - n_1(H)$. Further since four distinct edges are deleted from $H$ when constructing $H^*$, the contribution of these deleted edges to the sum $2m(H)$ is~$8$.

Let $z \in N_H(x) \setminus \{y^*\}$. Then, $d_{H^*}(z) \le
d_{H'}(z) \le 1$. If $d_{H^*}(z) = 1$, then $z$ contributes~$2$ to
the sum $2n(H) - n_1(H)$ and~$1$ to the sum $2n(H^*) - n_1(H^*)$. If
$d_{H^*}(z) = 0$, then $z$ contributes at least~$1$ to the sum
$2n(H) - n_1(H)$ and~$0$ to the sum $2n(H^*) - n_1(H^*)$ (since $z$
is deleted in~$H^*$). In both cases the contribution of $z$ to
$\Theta(H^*)$ is at least one less than its contribution to
$\Theta(H)$. This is true for every vertex in~$N_H(x) \setminus
\{y^*\}$. Hence the total contribution of the neighbors of $x$
different from~$y^*$ to $\Theta(H) - \Theta(H^*)$ is at
least~$|N_H(x) \setminus \{y^*\}| = |N_H(x)| - 1 \ge k \ge 5$.
Together with our earlier observation that the vertices $x$ and
$y^*$, together with the four edges incident with $x$ or $y^*$ in
$H$, contribute~$12$ to~$\Theta(H)$, this implies that $\Theta(H^*) \le \Theta(H) - 12 - 5 < \Theta(H) - 14 = \Theta(H) - 7|X^*|$, contradicting Claim~II.~\smallqed

\2
We now return to the proof of Theorem~\ref{t:kge5}. By Claim~IV, $\Delta(H) = 2$. Let $x$ be a vertex in $H$ with $d_H(x) = 2$. Since $\taut(H) \ge 3$ by Claim~I, and since $H$ is connected, there exists an edge, $e$, that contains a neighbor, $y$, of $x$ but does not contain~$x$. Let $X = \{x,y\}$ and consider the hypergraph $H - X$. Since $d_H(x) = 2$ and $d_H(y) = 2$, the vertices $x$ and $y$ both contribute~$2$ to the sum $2n(H) - n_1(H)$. Further the three edges incident with $x$ or $y$ contribute~$6$ to the sum $2m(H)$. Furthermore, each vertex in $N_H(x) \setminus \{y\}$ has degree~$0$ or~$1$ in $H - X$ and therefore contributes at least~$1$ to $\Theta(H) - \Theta(H - X)$. This implies that $\Theta(H - X) \le \Theta(H) - 10 - (|N_H(x)| - 1) \le \Theta(H) - 10 - k + 1 \le \Theta(H) - 14 = \Theta(H)  - 7|X|$, contradicting Claim~II. This completes the proof of Theorem~\ref{t:kge5}.~\qed

\2
As an immediate consequence of Theorem~\ref{t:kge5}, we have the following results.

\begin{cor}
For $k \ge 5$, if $H \in \cH_k$, then $7\taut(H) \le 2\nH + 2\mH$.
 \label{c:kge5}
\end{cor}

\begin{cor}
For all $k \ge 5$, we have $b_k \le \frac{2}{7}$.
 \label{c:bge5}
\end{cor}

Theorem~\ref{t:main2a} follows from Corollary~\ref{c:b2},  Corollary~\ref{c:b3b4} and Corollary~\ref{c:kge5}.

\subsection{Proof of Theorem~\ref{t:main1B}}
\label{S:main1B}

In this section, we present a proof of Theorem~\ref{t:main1B}. Recall its statement.
\2

\noindent \textbf{Theorem~\ref{t:main1B}}. \emph{For $k \ge 4$, if $H \in \cH_k^*$, then $\displaystyle{\gt(H) \le \left( \max \left\{ \frac{2}{k+2}, b_{k-1}\right \} \right) \nH .}$
}

%\newpage

\noindent \textbf{Proof of Theorem~\ref{t:main1B}.} Suppose to the contrary that the theorem is not true. Let $H \in \cH_k^*$ be a counterexample with $\nH + \mH$ a minimum. We proceed in a similar manner as in the proof of  Theorem~\ref{t:main1A}.

%\newpage
\begin{unnumbered}{Claim~I.}
The following properties hold in the hypergraph $H$. \\
\indent {\rm (a)} $H$ is connected. \\
\indent {\rm (b)} The deletion of any edge in $H$ creates an isolated vertex or an isolated edge. \\
\indent {\rm (c)} There is no dominating vertex in $H$.
\end{unnumbered}
\textbf{Proof of Claim~I.}  Parts~(a) and~(b) follows from the minimality of $H$ and the observation that the deletion of an edge cannot decrease the total domination number. To prove Part~(c), suppose that $H$ contains a dominating vertex $v$. The vertex $v$ and any one of its neighbors forms a TD-set in $H$, implying that $\gt(H) = 2$. By Part~(b), $H$ contains no isolated vertex or isolated edge. Since no two edges of $H$ intersect in $k-1$ vertices, we therefore have that $\nH \ge k+2$. Hence, $\gt(H) \le 2\nH/(k+2)$, contradicting the minimality of $H$. This proves Part~(c).~\smallqed

\begin{unnumbered}{Claim~II.}
Every edge in $H$ contains at least one degree-$1$ vertex.
\end{unnumbered}
\textbf{Proof of Claim~II.} We proceed as in the proof of Claim~\ref{claim2}. Let $u$, $v$, $e$, $e_1$ and $e_2$ be defined as in the proof of Claim~\ref{claim2}. If the edge $e_2$ contains a degree-$1$ vertex, then at least one vertex in addition to the vertices in $e_1 \cup \{u\}$ becomes an isolated vertex when we delete all edges incident with $u$ or $v$. Thus in this case we delete at least $k+2$ vertices and we add two vertices to $T$, and we proceed as in the 2nd paragraph of the proof of Claim~\ref{claim2}. In this case, $|T| = 2 + \ell$, where $\ell \ge 0$ denotes the number of isolated edges created when removing $u$ and $v$, and at least $k + 2 + k \ell$ vertices are deleted from $H$. Thus if $n'$ denotes the number of vertices in $H$ that are not deleted in the process, then $n' \le \nH - k - 2 - k \ell$, implying that

\[
\begin{array}{lcl} \3
\displaystyle{ \left( \frac{2}{k+2} \right) (\nH - n') } & \ge &
 \displaystyle{ \left( \frac{2}{k+2} \right) \left( k+2 + k \ell  \right) } \\ \3
& = & \displaystyle{  2 +  \left( \frac{2k}{k+2} \right) \ell }  \\ \3
& \ge & \displaystyle{  2 + \ell } \\ \3
& = & |T|.
\end{array}
\]

Suppose that the edge $e_2$ does not contain any degree-$1$ vertices. Then there is an edge, $e_3$, which would become isolated after the deletion of the edge $e_2$ from $H_2$. We note that neither $u$ nor $v$ belong to the edge $e_3$ and therefore that $e_3 \notin \{e,e_1,e_2\}$. Let $w \in e_2 \cap e_3$. We now delete all edges incident with a vertex in the set $\{u,v,w\}$ and we delete any resulting isolated vertices. Further we add the three vertices $u$, $v$ and $w$ to the set $T$. We note that every vertex in $e_1 \cup e_3 \cup \{u\}$ becomes an isolated vertex. We therefore delete at least $2k+1$ vertices and we add three vertices to $T$. If this process creates an isolated edge, then from each such isolated edge $f$, if any, we choose one vertex that is a neighbor of a vertex in $T$ and add it into $T$, and delete the $k$ vertices in $f$. Hence in this case, $|T| = 3 + \ell$, where $\ell \ge 0$ denotes the number of isolated edges created when removing $u$, $v$ and $w$, and at least $2k + 1 + k \ell$ vertices are deleted from $H$. Thus if $n'$ denotes the number of vertices in $H$ that are not deleted in the process, then $n' \le \nH - 2k - 1 - k \ell$. Since $k \ge 4$, we note that $2(2k+1)/(k+2) \ge 3$ and $2k/(k+2) > 1$, implying that

\[
\begin{array}{lcl} \3
\displaystyle{ \left( \frac{2}{k+2} \right) (\nH - n') } & \ge &
 \displaystyle{ \left( \frac{2}{k+2} \right) \left( 2k + 1 + k \ell \right) } \\ \3
 & = & \displaystyle{ \left( \frac{2(2k+1)}{k+2} \right)  + \left( \frac{2k}{k+2} \right) \ell   }  \\ \3
& \ge & \displaystyle{  3 + \ell } \\ \3
& = & |T|.
\end{array}
\]

In both cases, we therefore have that $|T| \le 2(\nH - n')/(k+2)$. If $n' = 0$, then the set $T$ is a TD-set in $H$, implying that $\gt(H) \le |T| \le 2\nH/(k+2)$, a contradiction. Hence, $n' > 0$. Let $H'$ denote the resulting hypergraph on these $n'$ vertices. Let $H'$ have size~$m'$. By construction, the hypergraph $H'$ is in the family~$\cH_k^*$. In particular, we note that $n' \ge k+2$. By the minimality of $H$, we have that
\[
\gt(H') \le \left( \max \left\{ \frac{2}{k+2}, b_{k-1}\right \} \right) n'.
\]

Let $T'$ be a $\gt(H')$-set and note that the set $T \cup T'$ is a TD-set of $H$. Suppose that $2/(k+2) \ge b_{k-1}$. Then, $|T'| \le 2n'/(k+2)$, and so

\[
\gt(H) \le |T \cup T'| \le \left( \frac{2}{k+2} \right) (\nH - n') + \left( \frac{2}{k+2} \right) n' = \left( \frac{2}{k+2} \right) \nH,
\]

a contradiction. Hence, $2/(k+2) < b_{k-1}$. Thus, $|T'| \le b_{k-1}n'$, and so
\[
\gt(H) \le |T \cup T'| \le \left( \frac{2}{k+2} \right) (\nH - n') + b_{k-1}n' < b_{k-1} (\nH - n') + b_{k-1}n'  = b_{k-1} \nH,
\]

a contradiction. This completes the proof of Claim~II.~\smallqed

\2
We now return to the proof of Theorem~\ref{t:main1B}. By Claim~II, every edge in $H$ contains at least one degree-$1$ vertex. Let $H'$ be the hypergraph obtained from $H$ by deleting exactly one degree-1 vertex from each edge. Since $H \in \cH_k^*$, we note that no multiple edges are created. Further, $H'$ contains no isolated edge and no isolated vertices, and so $H' \in \cH_k^*$. Let $H'$ have order $n'$ and size $m'$. Then, $n' = \nH - \mH$ and $m' = \mH$. Thus, $n' + m' = \nH$. We note that $k - 1 \ge 3$. By Definition~\ref{defn2} we have that $\taut(H') \le (n' + m') b_{k-1} \le b_{k-1} \nH$. An identical proof as in the proof of Claim~\ref{claim4} of Theorem~\ref{t:main1A} shows that $\gt(H) = \taut(H')$, implying that $\gt(H) \le b_{k-1} \nH$, a contradiction. This  completes the proof of Theorem~\ref{t:main1B}.~\qed

\section{Tight Asymptotic Bounds}
\label{S:asym}

In this section we prove Theorem~\ref{t:asym} which establishes a tight asymptotic upper bound on $b_k$ for $k$ sufficiently large.
Since every strong transversal in a hypergraph, $H$, is a total transversal in $H$, and since every total transversal in $H$ is a transversal in $H$, we have the following observation.

\begin{ob}
For every hypergraph $H$, we have $\tau(H) \le \taut(H) \le \ST{H}$. \label{ob:chain}
\end{ob}

Using probabilistic arguments, Alon~\cite{Al90} established the following result.

\begin{thm}{\rm (\cite{Al90})}
For every $\varepsilon > 0$ and sufficiently large $k$ there exist $k$-uniform hypergraphs, $H$, satisfying
\[
\tau(H) \ge \left(\frac{ (1-\varepsilon) \, \ln(k)}{k} \right)(\nH +
\mH)
\]
 \label{th:Alon}
\end{thm}

The following result establishes a tight asymptotic upper bound on the strong transversal number of a $k$-uniform hypergraph for $k$ sufficiently large.

\begin{thm}
 For every constant $c> 1$ and every $k$-uniform hypergraph $H$, we have
\[
\ST{H} \le \left( \frac{ \ln(k) + \ln(c) }{k-1} \right) \nH + \left( \frac{ \ln(k) + \ln(c) }{c(k-1)} \right) \mH   + \left( \frac{ 2 }{ck} \right) \mH .
\]
 \label{th:strong}

\end{thm}
\proof  Let $H = (V,E)$
and let $p = \ln(ck) / (k-1)$. Let $X_1$ be a random subset of
$V(H)$ where a vertex {$x$ is chosen to be in $X_1$ with probability $\Pr(x\in X_1) = p$, independently of the choice for any other vertex. For every edge $e \in E$ that does not intersect
$X_1$, select two vertices from $e$ and let $X_2 \subseteq V$ be the
resulting set of all such selected vertices. For every edge $e \in
E$ such that $|e \cap X_1| = 1$, select one vertex from $e \setminus
X_1$ and let $X_3 \subseteq V$ be the resulting set of all such
selected vertices. The resulting set $X_1 \cup X_2 \cup X_3$ is a
strong transversal in $H$. The expected value of the set $X_1$ is
\[
\bbE(|X_1|) = p \nH = \left( \frac{ \ln(k) + \ln(c) }{k-1} \right) \nH.
\]

Using the inequality $1 - x \le {\rm e}^{-x}$ for $x \in \mathbb{R}$, the expected value of the set $X_2$ is given by

\[
\begin{array}{lcl}
\bbE(|X_2|) & \le & (1 - p)^k  \cdot \mH  \cdot 2 \2 \\
& = & \displaystyle{ \left( 1 - \frac{ \ln(ck) }{k-1} \right)^k \cdot 2 \mH  } \2 \\
& = & \displaystyle{ \left( \left( 1 - \frac{ \ln(ck) }{k-1} \right)^{\frac{k-1}{\ln(ck)}} \right)^{\frac{k}{k-1} \ln(ck)} \cdot 2 \mH  } \2 \\
& < & \displaystyle{ {\rm e}^{-\frac{k}{k-1} \ln(ck)} \cdot 2 \mH  } \2 \\
& \le & \displaystyle{ \frac{2}{ck} \cdot \mH  }.
\end{array}
\]

The expected value of the set $X_3$ is given by

\[
\begin{array}{lcl}
\bbE(|X_3|) & \le & \mH \cdot k \cdot p \cdot (1 - p)^{k-1} \2 \\
& = & \displaystyle{ k \left( \frac{ \ln(ck) }{k-1} \right) \left( 1 - \frac{ \ln(ck) }{k-1} \right)^{k-1} \cdot \mH }
\2 \\
& = & \displaystyle{ k \left( \frac{ \ln(ck) }{k-1} \right) \left( \left( 1 - \frac{ \ln(ck) }{k-1} \right)^{\frac{k-1}{\ln(ck)}} \right)^{\ln(ck)} \cdot \mH  } \2 \\
& < & \displaystyle{ k \left( \frac{ \ln(ck) }{k-1} \right) {\rm e}^{-\ln(ck)} \cdot \mH  } \2 \\
& = & \displaystyle{ \left( \frac{ \ln(k) + \ln(c) }{c(k-1)} \right) \mH  }. \\
\end{array}
\]

By linearity of expectation, we have that
 $\bbE(|X_1 \cup X_2 \cup X_3|) \le \bbE(|X_1|) + \bbE(|X_2|) + \bbE(|X_3|)$,
  yielding the desired upper bound.~\qed

\2
As a consequence of Theorem~\ref{th:strong}, we have the following results.

\begin{cor}
Given any $\varepsilon > 0$, if $H$ is a $k$-uniform hypergraph with $k$ sufficiently large, then
\[
\ST{H} < \left( (1+\varepsilon) \, \frac{ \ln(k)}{k} \right) (\nH + \mH).
\]
 \label{c:strong}
\end{cor}
\proof For a constant $c > 1$, we note that the functions,
\[
\frac{ \ln(k) + \ln(c) }{k-1} \hspace*{0.5cm} \mbox{and}
 \hspace*{0.5cm}  \frac{ \ln(k) + \ln(c) }{c(k-1)} +\frac{2}{ck},
\]
tend to $\ln(k)/(k-1)$ and $\ln(k)/(c(k-1)) < \ln(k)/(k-1)$, respectively, when $k$ tends to infinity.
 Hence for $k$ sufficiently large, we have that
\[
\max \left\{ \frac{ \ln(k) + \ln(c) }{k-1},   \frac{ \ln(k) + \ln(c) }{c(k-1)} + \frac{2}{ck}  \right\} < (1+\varepsilon) \, \frac{ \ln(k)}{k}.
\]

\noindent
Therefore for $k$ sufficiently large, we have that

\[
\left( \frac{ \ln(k) + \ln(c) }{k-1} \right) \nH + \left( \frac{ \ln(k) + \ln(c) }{c(k-1)} \right) \mH   + \left( \frac{ 2 }{ck} \right) \mH < \left( (1+\varepsilon) \, \frac{ \ln(k)}{k} \right) (\nH + \mH).
\]

\noindent
The desired result now follows from Theorem~\ref{th:strong}.~\qed

\2
We are now in a position to prove Theorem~\ref{t:asym}. Recall its statement.

\medskip
\noindent \textbf{Theorem~\ref{t:asym}}. \emph{For $k$ sufficiently large, we have that
$\displaystyle{ b_k = (1 + \oo(1)) \frac{ \ln(k)}{k}}.$ }

\2
\noindent \textbf{Proof of Theorem~\ref{t:asym}.} It suffices for us to prove that for $\varepsilon > 0$ and for $k$ sufficiently large, we have
\[
(1-\varepsilon) \, \frac{ \ln(k)}{k}  \le b_k \le  (1+\varepsilon) \, \frac{ \ln(k)}{k}.
\]

The upper bound on $b_k$ follows from Observation~\ref{ob:chain} and Corollary~\ref{c:strong}.
For the lower bound let $\varepsilon>0$ and let $k$ be sufficiently large, such that a $k$-uniform hypergraph, $H$, exists with
$\tau(H) \ge [(1-\varepsilon) \, \ln(k)/k] (\nH + \mH)$ (which exists by Theorem~\ref{th:Alon}).
 Assume that $H$ contains $n_0$ isolated vertices and $e_0$ isolated edges. Let $H'$ be obtained from $H$ by deleting all isolated vertices and isolated edges and the vertices belonging to isolated edges. Then, $H' \in \cH_k$. Further, $n(H') = \nH - n_0 - ke_0$ and $m(H') = \mH - e_0$. As $n_0 \ge 0$ and $(1-\varepsilon) \ln(k) (k+1) / k > 1$ when $k$ is sufficiently large, we have that

\[
\begin{array}{lcl} \2
\taut(H')
& \ge & \tau(H')  \\ \2
& = & \tau(H) - e_0 \\ \3
& \ge & \displaystyle{ \left( \frac{(1-\varepsilon) \, \ln(k)}{k} \right) (\nH + \mH) - e_0 } \\ \3
& \ge &
 \displaystyle{ \left( \frac{(1-\varepsilon) \, \ln(k)}{k} \right) (n(H') + m(H') + n_0 + ke_0 + e_0 ) - e_0 } \\ \2
 & = & \displaystyle{ \left( \frac{(1-\varepsilon) \, \ln(k)}{k} \right) (n(H') + m(H')) } \\ \3
  & & \hspace*{0.5cm} \displaystyle{ + \left( \frac{(1-\varepsilon) \, \ln(k)}{k} \right) (n_0 + (k+1)e_0) - e_0. }  \\ \3
& \ge & \displaystyle{ \left( \frac{(1-\varepsilon) \, \ln(k)}{k} \right) (n(H') + m(H')) }. \\
\end{array}
\]

\noindent
This implies that $b_k \ge \left( \frac{ (1-\varepsilon) \, \ln(k)}{k} \right)$, which
establishes the desired lower bound on~$b_k$ and completes the proof of Theorem~\ref{t:asym}.~\qed

\section{Closing Remarks and Open Problem}

In view of Theorem~\ref{t:main1A}, it is of interest to determine the value of $b_k$ for $k \ge 2$. In Theorem~\ref{t:main2a} we show that $b_2 = \frac{2}{5}$ and $b_3 = \frac{1}{3}$, and we show that $b_{k-1} \le 2/(k+1)$ for $k \in \{3,4,5,6\}$.  In Theorem~\ref{t:asym}, we establish a tight asymptotic bound on $b_k$ for $k$ sufficiently large which shows that is not true that $b_{k-1} \le 2/(k+1)$ when $k$ is large enough. We pose the following problems that
still remain to be settled.

\begin{prob}
Determine the exact value of $b_k$ for $k \ge 4$. \label{P1}
\end{prob}

\begin{prob}
Determine the smallest value of $k$ for which $b_{k-1} > 2/(k+1)$.
 \label{P3}
\end{prob}

\end{document}